\documentclass[12pt]{amsart}

\usepackage{amsmath,amsxtra} 
\usepackage{amssymb,amsthm,amsfonts}
\usepackage{latexsym}
\usepackage{graphicx}
\usepackage{epsfig}

\newtheorem{theorem}{Theorem}

\newcommand{\cal}{\mathcal}

\newcommand{\til}{\tilde}

\newcommand{\n}{\newline}
\newcommand{\nd}{\noindent}

\newcommand{\R}{\mathbb{R}}
\newcommand{\C}{\mathbb{C}}

\def\vphi{\varphi}
\def\lam{\lambda}

\begin{document}

\title[1D Schr\"odinger Operators]%
      {A representation formula related to
   Schr\"odinger operators} 
\author[S.~Zheng]%
       {Shijun Zheng}
\address{Department of Mathematics \\
        Louisiana State University \\
        Baton Rouge, LA 70803}
\email{szheng@math.lsu.edu}
 \urladdr{http://www.math.lsu.edu/\textasciitilde{szheng}}
 \thanks{The author is supported in part by DARPA 
(Defense Advanced Research Projects Agency)}         
\date{\today}

\begin{abstract}

Let $H=-d^2/dx^2 +V$ be a Schr\"odinger operator on the real line, 
where $V\in L^1 \cap  L^2$.  We define the perturbed Fourier transform 
$\mathcal{F}$ for $H$ and show
that $\mathcal{F}$ is an isometry from the absolute continuous
subspace  onto $L^2(\mathbb{R})$. 
This property allows
us to construct a kernel formula for the spectral operator $\vphi(H)$.  


\end{abstract}

\maketitle 


Schr\"odinger operator is a central subject in the mathematical study
of quantum mechanics. 
Consider the Schr\"odinger operator $H=-\triangle +V$ on $\mathbb{R}$,
where $\triangle=d^2/dx^2$ and the potential function $V$ is real valued.
In Fourier analysis, it is well-known that a square integrable function admits 
an expansion with exponentials as eigenfunctions of $-\triangle$. A
natural conjecture is that an $L^2$ function admits a similar 
expansion in terms of
``eigenfunctions'' of $H$, a perturbation of the Laplacian (see \cite{RS}.  Ch.XI and the notes), under certain condition on $V$.

The three dimension analogue was proven true by T.Ikebe \cite{Ik}, a
member of Kato's school, in 1960.  Later his result was extended by
Thor to
the higher dimension case \cite{Th}. 
In one dimension, recent related results
can be found in e.g., Guerin-Holschneider \cite{GH}, Christ-Kiselev \cite{ChK}
and Benedetto-Zheng \cite{BZhen}.

Throughout this paper we assume $V: \R\rightarrow \R $ is in $L^1\cap 
L^2$. We shall prove a one-dimensional version of Ikebe's
theorem for $L^2$ functions (Theorem 1).  
Theorem 2 presents an integral formula for the kernel 
of the spectral operator 
$\vphi(H)$ for a continuous function $\vphi$ with compact 
support. In a sequel to this
paper we shall use this explicit formula to study function spaces
associated with $H$ (see  \cite{BZhen}).





The {\em generalized eigenfunctions} $e(x,\xi),\; \xi\in \R$ of $H$ satisfy
\begin{equation}
(-d^{\,2}/dx^2 + V(x)) e(x,\xi) = \xi^2 e(x,\xi) 
\label{S_eq} 
\end{equation} 
in the sense of distributions. 

\noindent
Definition.
  The {\em perturbed Fourier 
transform} $\mathcal{F}$ on $L^2$ is given by  
\begin{align}
\mathcal{F}f(\xi) = \textrm{l.i.m}.(2\pi)^{-1/2}\int f(x) 
\overline{e(x,\xi)} \, dx\\
=\lim_{N\rightarrow \infty} (2\pi)^{-1/2}\int_{-N}^N f(x) 
\overline{e(x,\xi)} \, dx,  \notag 
\end{align}
where the convergence is in $L^2$ norm as $N \rightarrow \infty$.  
By Theorem 1, $\cal{F}$ is a well-defined isometry from $\cal{H}_{ac}$ onto 
$L^2$.

\begin{theorem} Suppose
$V\in L^1\cap L^2$.  Then there
exists a family of solutions $e(x,\xi)$, 
$|\xi| \in [0,\infty)\setminus\cal{E}_0$, $\cal{E}_0 $ being a bounded closed set of measure zero,
to equation   
(\ref{S_eq}) with the following properties.  \n\n
(i) If $f \in L^2$, then 
 there exists an element $\til{f} \in L^2$ such that
$$ \cal{F}f(\xi)= \tilde{f} (\xi)
\qquad \textrm{in}\; L^2.  $$  
(ii) The adjoint operator $\cal{F}^*$ is given by 
$$
\cal{F}^*g=\textrm{l.i.m.}_{N\rightarrow\infty}\; \sum_{i=1}^N (2\pi)^{-1/2} \int_{\alpha_i \leq \xi^2 \leq \beta_i}
g(\xi)e(x,\xi) \, d\xi 
$$
in $L^2$,
 where  $[\alpha_i,\beta_i) \subset (0,\infty)$ are
a countable collection of disjoint intervals 
with $[0,\infty)\setminus\cal{E}_0^2$ equal to $\cup_i [\alpha_i,\beta_i)$.\n\n 
(iii) If $f \in L^2$, then $\Vert P_{ac}f \Vert_{L^2} = \Vert \tilde{f}
\Vert_{L^2}$, where $P_{ac}$ is the projection onto $\cal{H}_{ac}$, the  
absolute continuous subspace in $L^2$.   \n\n
(iv) $\cal{F}: L^2\rightarrow L^2$ is a surjection.
Moreover, $\cal{F}\cal{F}^*=Id$ and  $\cal{F}^*\cal{F}=P_{ac}$.\n\n  
(v) If $f \in \mathcal{ D}(H)$, then ${(Hf)}\sptilde(\xi) = \xi^2 \til{f}(\xi)$
in $L^2$. 
\end{theorem} 

\nd
Remark 1. 
The proof is based on the ideas of \cite{Ik} for $3D$.  We also use some
simplifications as found in Reed and Simon(\cite{RS}) and 
Simon \cite{Si2}. 
 
\nd
Remark 2. If $ \vert e(x,\xi)\vert \le C $ a.e. $(x,\xi)\in \R^2$, then
we have a ``better-looking'' form in ($ii$) of the theorem 
$$ \cal{F}^*g = \textrm{l.i.m}.(2\pi)^{-1/2} \int g(\xi)e(x,\xi) \, d\xi. 
$$








\vspace{.3in}

If $H=\int \lam dE_\lam$ is the spectral resolution of $H$, 
define the spectral operator $\vphi(H):=\int \vphi(\lam) dE_\lam $ 
by  functional calculus.
We prove a representation formula for  the integral
kernel of $\varphi(H)$. 




Let $\{e_k\}_{k=1}^\infty$ be an orthonormal basis in $\cal{H}_p$, the 
subspace of eigenfunctions in $L^2$ for $H$ and let $\lam_k$ be the eigenvalue corresponding to $e_k$.

\begin{theorem}    
Let the operator $H$ be as in Theorem 1. 
Suppose $\varphi : \R \rightarrow \C$ is continuous  and has a compact 
 support disjoint from $\cal{E}_0^2:=\{\eta^2:\eta\in \cal{E}_0  \}$. 
Then for $f \in L^1 \cap L^2$
\begin{equation}
\varphi(H)f(x) = \int_{-\infty}^\infty K(x,y) f(y) \, dy 
\label{eq:kernel}
\end{equation}
where
$K=K_{ac}+K_p$,
$$
K_{ac}(x,y) = (2\pi)^{-1} \int_{-\infty}^\infty \varphi(\xi^2) e(x,\xi)
\overline{e(y,\xi)} \, d\xi.
$$
 and 
$$ 
K_p(x,y)= \sum_k \vphi(\lam_k)e_k(x)\bar{e}_k(y).
$$
\end{theorem}

\nd
Remark 1. If $|e(x, \xi) | \le C$, a.e. $(x,\xi) \in \R^2$, 
then, under the same condition the integral expression (\ref{eq:kernel}) is
valid for any $\vphi \in C(\R)$ with compact support. 

\nd
Remark 2. When $\vphi$ is smooth with rapid decay and $V$ is compactly
supported in $\R^3$, a formula
of this type appeared in \cite{Tao} by Tao.

\end{document}